\documentclass[12pt]{article}
\usepackage{amsfonts}

\begin{document}
\centerline{ A new inverse formula for the Laplas`s
transformation. }

\centerline{ Pavlov An.V.(MIREA(TU)).}

 In the article is proved,that the complex part of the analytical
continuation of the  $$r(p)=\mathcal L \mathcal L Z(x)=
\int\limits_{0}^{\infty}e^{-pt}dt\int\limits_{0}^{\infty}e^{-tx}Z(x)
 dx, p\in\{p:Im\,p\geq0\}, $$
               equals to $
-\pi Z(x),x\in(0,\infty),$ if $p=s=-x\in(-\infty,0)$ for a wide
class of a functions $Z(x):$

It is proved,that  the odd functions $$Z(x)=\sum_{k=1}^{l}\gamma_k
e^{\lambda_k x},\,\,\gamma_k=
res_{p=\lambda_k}\frac{Q_n(p)}{P_l(p)},\,\lambda_k=-\alpha_k+\beta_k
i,\alpha_k>0,k =1,\ldots,l,$$  satisfy to all the conditions of
the theorems of the article.

THE TRANSFORMS OF FOURIER, THE TRANSFORM OF LAPLAS,THE NEW FORMULA
OF  TRANSFORMATION \\[1cm]

\bf Introduction. \rm

In the article we prove (a theorem 2),that for the negative
variable $p=-x,x\in(0,\infty)$ a complex part of the analytical
%  prolongation
continuation
              $r(p)=\mathcal L\mathcal L Z(x)$
               equals to $
-\pi Z(x),x\in(0,\infty)$ for a wide class of a functions $Z(x):$

$$\pi Z(s)=-Im\,r_{\mathcal An}(-s),\,s\in(0,\infty)\,\,
,r_{\mathcal An}(s)\equiv r(s),s\in(0,+\infty),$$ $$ r(s)=\mathcal
L \mathcal L (Z(x)),\,s\in(0,\infty),\,\mathcal L(Z(x))
=\int\limits_{0}^{\infty}e^{-sx}Z(x)dx,\, s\in(0,\infty), $$
where,by definition,$r_{\mathcal An}(p),p\in\mathbb C,$  is the
analytical continuation of the function $r(p),p\in
D=\{p:Re\,p>0\},$ and the function $r_{\mathcal An}(p)$ is regular
(analytic) in $\mathbb C\setminus P,\,\,P=\{z_j:z_j\in \mathbb
C,z_j \notin(-\infty,0) ,j=1,\ldots,m\},$ the set P may be
$\emptyset.$

All results of the article are true, for instance,if
$$Z(x)=\sum_{k=1}^{2N}\gamma_k e^{\lambda_k
x}=\sum_{k=1}^{N}\gamma_k e^{-\alpha x}\cos {\beta_k
x}=Re\,Z(x),\,\lambda_k=-\alpha_k+\beta_k
i,\,\,\alpha_k\in(0,\infty),\,\gamma_k,\beta_k\in(-\infty,\infty),$$
$\lambda_{2j+1}=-\alpha_{2j+1}+\beta_{2j+1}i,\lambda_{2j+2}=
-\alpha_{2j+1}-\beta_{2j+1}i,\,\,j=1,\ldots,N-1,\,\,\lambda_j\neq
\lambda_i, $ for all $i\neq j$.
%Этому классу
%принадлежат также функции такого вида,умноженные на произвольную
%непрерывную на $[0,\infty)$ функцию $z(x): z(x)\rightarrow0,$  при
%$x\rightarrow0.$

The equality  $r(s)=-i\mathcal F_0\mathcal
 F_0(Z(x))=\mathcal L\mathcal L(Z(x)),s\in(0,+\infty),$ is considered
 in the theorem 1,where,by definition,
$$ \mathcal
F_0(Z(x)))=\int\limits_{0}^{\infty}e^{iyx}Z(x)dx,\,\,y\in(-\infty,\infty).
$$ The main result of the article (the theorem 2) follows from the
theorem 1.

 All the conditions of the theorems 1 and 2 are checking
in the remark 1 for the above functions $Z(x).$

\bf 1. The main result. \rm

We shall use a designations  $f(p)\in\mathcal An\,G,$ if the
function  $f(p)$ is regular in the open domain      $G\in\mathbb
C.$

 Theorem 1.

If for the function $Z(p),$ we have$$R(p)=
\int\limits_{0}^{\infty}e^{-px}dx\int\limits_{0}^{\infty}e^{itx}Z(t)dt\in
\mathcal An\, D_{*}=\{p:-\pi/2< arg \,p<\varphi_0 \}, \, $$
$$F(p)=\int\limits_{0}^{\infty}\frac{e^{-x_1}}{p}dx_1
\int\limits_{0}^{\infty}e^{itx_1/p} Z(t)dt\in \mathcal An
\,D_{*},$$ for a variable $\varphi_0\,\,:0<\varphi_0<\pi/2;$ if
for all $p=-is,s\in(0,\infty),$ the functions are continuous.

Then $$\mathcal L\mathcal L(Z(x))=-i\mathcal F_0\mathcal
F_0(Z(x))=
Im\,\int\limits_{0}^{\infty}e^{isx}dx\int\limits_{0}^{\infty}e^{itx}Z(t)dt,
\,\,s\in(0,\infty).$$

Proof.

We can replace  $ux=x_1,p=s\in(0,\infty),$ in the integral of
$R(p),$ if $p=u\in(0,\infty).$  We obtain
$R(u)=F(u),s\in(0,\infty).$ The functions are regular in
$D_{*}\{p:-\pi/2< arg \,p<\varphi_0 \}, \, $ it is the open domain
and $(0,\infty)\in D_{*},$ then  $R(p)=F(p),p\in D.$ The function
are continuous
%on a boundary
for $p=is,s\in(-\infty,0),$
                                               then
 $R(-is)=F(-is),s\in(0,\infty),$ and after the inverse replace
 $x_1/s=x,s\in(0,\infty),x_1\in[0,\infty)$ in the integral of $F(-is)$
we obtain the theorem 1..

Theorem 2.

Let the function $\mathcal F_0\mathcal F_0(Z(x))=r_{F}(p)$ is
regular in the domain $D_{+}=\{p:Im\,p>0\}$ and is continuous for
all real  $p=s\in(-\infty,\infty).$

Let all the conditions of the theorem 1 are holds.

Let the function $r(s)=\mathcal L \mathcal
L(Z(x)),s\in(0,\infty),$ can be analytical continued in the left
part of the plane : $
D_{-}=\{p:Re\,p<0\setminus\{p_j,j=1,\ldots,m\}\,\,\}$ ( $m$ can be
0).

then $$-\pi Z(x)= Im\,r_{\mathcal
An}(-x),x\in(0,\infty),\,\,r_{\mathcal An}(p)=r(p),p\in
D=\{p:Re\,p>0\}. $$

Proof.

Let

 $$-ir_{F}(p)=-i\mathcal F_0\mathcal F_0(Z(x)),
r_{L}(p)=r_{\mathcal An}(p),r_{L}(p)= \mathcal L\mathcal
L(Z(x)),p=s\in(0,\infty).$$

We can see,that the real part of $r_{F}(s),s\in(-\infty,\infty),$
is equal to $$Re\,r_{F}(s)=
\int\limits_{0}^{\infty}\cos{sx}\,dx\int\limits_{0}^{\infty}\cos{tx}\,
Z(t)dt-\int\limits_{0}^{\infty}\sin{sx}\,dx\int\limits_{0}^{\infty}
\sin{tx}\, Z(t)dt=$$ $$
=(1/4)\int\limits_{-\infty}^{\infty}e^{isx}\,dx
\int\limits_{-\infty}^{\infty}e^{itx} \,
[Z_{od}(t)+Z_{ev}(t)]dt=(1/4)\int\limits_{-\infty}^{\infty}e^{isx}\,dx
\int\limits_{-\infty}^{\infty}e^{itx} \, Z_{+}(t)dt=(\pi/2)
Z_{+}(-s),$$ $s\in(-\infty,\infty),$ where $$Z_{od}(x)\equiv
Z_{ev}(x)\equiv Z_{+}(x)/2=Z(x),x\in(0,+\infty);$$
$$Z_{od}(-x)\equiv
Z_{od}(x),Z_{ev}(-x)\equiv-Z(x),Z_{+}(-x)\equiv0,-x\in(-\infty,0).$$

%We have chack

We obtain   $$ Im\,(-ir_{F}(-s))\equiv
 (-\pi/2)Z_{+}(-(-s))=-\pi Z(s),\,s\in(0,\infty)
 \eqno(1.1).$$

From the theorem 1 we have $$r_{L}(s)=
-ir_{F}(s)=r(s),s\in(0,+\infty).$$

We shall prove,that the functions are equal for all the complex
 $p:Im\,p\geq0.$

 Let $p=z_j$ is a point,where the $r_L(p)$  or the $r_{F}(p)$
 functions
 are not regular; let $Q=\{p=z_j,j=1,\ldots,M\}$ is the set of all the
 points.

All the values of the function $ir_{F}(p)$ are real,if
$p=s\in(0,\infty)$ -  it  follows from the the theorem 1,where
%the aprovement of the theorem 1
 $$ -ir_{F}(p)\equiv Re\,(-ir_{F}(p))\equiv \mathcal L\mathcal
L(Z(x))=r(s),p=s\in(0,\infty).$$ Then,we can use the Reaman`s
theorem about the symmetrical continuation ([2]),and the function
$-ir_{F}(p)= -i\mathcal F_0\mathcal F_0(Z(x)),Im\,-ir_{F}(s)\equiv
0,s\in(0,\infty),$ can be analytical continued from the domain
$D_0=\{p:Im\,p>0\}$ in the the low part  of the complex plane in
the domain $D_1,$ such that the real axis $(0,\infty)\in D_1:$ :
$(0,\infty)\in D_1,-ir_{\mathcal An}(p)\in\mathcal An D_1,D_0\in
D_1,-ir_{\mathcal An}\equiv -ir_{F}(p),p\in D_0$.

Now,from $-ir_{\mathcal An}(p)=r_{L}(p),p=s\in(0,\infty)\in D_1,$
it follows $-ir_{\mathcal An}(p)=r_{L}(p),p\in D_1\setminus Q.$

As the result we have $-ir_{\mathcal An}(p)=-ir_{F}(p),p\in
D_0\setminus Q\in D_1.$ (We use,that $-ir_{F}(p)\in \mathcal An
D_0\setminus Q.$

Therefore $r_{L}(p)=-ir_{\mathcal An}(p)=-ir_{F}(p),p\in
D_0\setminus Q.$

The function $r_{F}(p)$ is continuous ,if
$p=s\in(-\infty,\infty).$ Then, $$-ir_{F}(-s)=\lim_{p\rightarrow
-s}(-ir_{F}(p))=\lim_{p\rightarrow-s}r_{L}(p)=r_{L}(-s),\,s\in(0,\infty).$$

From the equality (1.1) we obtain $$-\pi
Z(s)=Im\,(-ir_{F}(-s))=Im\,r_{L}(-s)=Im\,r_{\mathcal
An}(-s),s\in(0,\infty).$$

The theorem 2 is proved.

Remark 1.

For all $l>n+2>1,l,n\in 1,2\ldots,$ ,we have (it is  a wallknown
formula)  $$\frac{Q_n(p)}{P_l(p)}=\frac{[q_0+q_1p+\ldots+q_n
p^n}{\prod_{k=1}^{l} (p-\lambda_k)}= 2\pi i
\int_{0}^{\infty}[\sum_{k=1}^{l}res_{p=\lambda_k}e^{p
x}\frac{Q_n(p)}{P_l(p)}]e^{-px}dx,\,Re\,p\geq0,$$ if
$\lambda_k\neq0,k=1,\ldots,l.$

  For all  complex $\lambda_{k}=-\alpha_{k}+\beta_k
i,\,\,\alpha_m\in(0,\infty),\,\beta_k\in(-\infty,\infty),$ the
function
% $$\frac{Q_n(t)}{\prod_{k=1}^{2N}(t-\lambda_k)}=
%\int_{0}^{\infty}[\sum_{k=1}^{2N}res_{p=\lambda_k}
$$Z(x)=\sum_{k=1}^{l}\gamma_k e^{\lambda_k x},\,\,\gamma_k=
res_{p=\lambda_k}\frac{Q_n(p)}{P_l(p)},\,k=1,\ldots,l$$ satisfy to
all the conditions of the theorems 1 and 2,if,for instance,
$\lambda_i\neq \lambda_j$  for all $i\neq j$.

(In the conditions of introduction  the function $Z(x)$ is equal
to $Z_{re}(x)$ ).

 Proof.

1. It is obviously,the function  $R(p),p=x+iy,$ is regular for all
 $x>0,y\in(-\infty,\infty)$ and continuous for all
$p=iy,y\in(-\infty,\infty)$
$$R(p)=\int\limits_{0}^{\infty}\frac{Q_n(t)}{P_l(t)}e^{-(x+iy)t}dt.$$

 2.

 If $p\in D_{*}=\{p=x+iy:-\pi/2<arg p
 <\varphi_0>0\},$ and
$0<\varphi_0<\min_{0<k<l}||\arctan \alpha_k/\beta_k|,$

 then
$$\frac{d\int\limits_{0}^{\infty}e^{-x_1}
 \frac{Q_n(-ix_1/(x+iy))}{\prod_{k=1}^{l}
((-ix_1/(x+iy))+\alpha_k-\beta_ki)}dx_1}{dp}=$$ $$=
\int\limits_{0}^{\infty}\frac{de^{-x_1}
 \frac{Q_n(-ix_1/p}{\prod_{k=1}^{l}
((-ix_1/p)+\alpha_k-\beta_ki)}}{dp}dx_1$$ is continuous for all
$p\in D_{*};$ we use $$
[x_1/\sqrt{x^2+y^2}](-i)(x-iy)=T(-ix-y)\neq-\alpha_k+\beta_k
i,\,\mbox{ for all } T\in(0,\infty), $$  if ,other $p\in
D_1=\{x+iy:x>0,y\leq 0\}$,( $\alpha>0),$ or $p\in
D_2=\{x+iy:y>0,x>0\bigcap-\pi/2-arg(-\alpha_k-|\beta_k|i)
<arg(-ix-y)<-\pi/2,k=1,\ldots,l\},$ where
$D_2=\{x+iy:x>0,y>0\bigcap
 0<arg (x+iy)<\min_{l=1,\ldots,l}(\arctan|\alpha_k/\beta_k|\},$
где $D_1\bigcup D_2=D_{*}.$

Then $$F(p)=(1/p)\int\limits_{0}^{\infty}e^{-x_1}
 \frac{Q_n(-ix_1/(x+iy))}{\prod_{k=1}^{l}
((-ix_1/(x+iy))+\alpha_k-\beta_ki)}dx_1\in\mathcal An D,$$

3. To prove,that $r_{\mathcal Am(p)} $ is continuous for all
$p=x\in(-\infty,0)$ and $r_{\mathcal Am}(p)=r(p)$,if $p\in
D_3=\{p:Im\,p>0\},$ we consider an equality  $$r(p)=\mathcal
L\mathcal
L(Z(x))=\int\limits_{0}^{\infty}Z(x)[1/(p+x)]dx\in\mathcal
An\mathbb C\setminus(-\infty,0).$$

We can write$$r(p)=\mathcal L\mathcal
L(Z(x))=\int\limits_{0}^{\infty}Z(x)[1/(p+x)]dx\in \mathcal An
\mathbb C\setminus(-\infty,0).$$

Let $$p=p_0 t,t\in(0,\infty),p_0=const.,arg\,
p_0=\pi-\varepsilon,0<\varepsilon,\ll 1,$$ and
 $x/(t p_0)=z,x\in[0,\infty),t\in(0,\infty).$
We have
$$r_1(p_0t)=\int\limits_{l(p_0)}Z(p_0tz)[1/(1+z)]dz,t\in(0,\infty),
\,l(p_0)=\{z:z=\tau/ p_0,\tau\in[0,\infty)\},$$ $(l(p_0)=\{z:arg\,
z=arg(1/p_0)\}),$ where
 $r_1(p)=r(p),p=p_0t,t\in(0,\infty);$  and
$r_1(p)\in\mathcal An D_{\varepsilon}=\{p:-2\varepsilon<arg\,
p/p_0<2\varepsilon\bigcap p\neq 0\},$
 while
$$r_1(p)=\int\limits_{l(p_0)}Z(pz)[1/(1+z)]dz=
 \int\limits_{0}^{\infty}Z((p/p_0)\tau)[1/(1+(\tau/p_0))]d\tau/p_0,$$
and,if $p\in D_{\varepsilon},$ we obtain $$\frac{dr_1(p)}{dp}=
 \int\limits_{0}^{\infty}\frac{dZ((p/p_0)\tau)}{dp}
 [1/(1+(\tau/p_0))]d\tau/p_0,$$ is continuous for all $p\in D_{\varepsilon},$
 $(Re\,(p/p_0)>const.>0).$

We obtain,that $r_1(p)=r(p),p\in l(p_0)),$  and
$r_1(p)=r_{\mathcal An}(p)=r(p) ,p\in
D_{\varepsilon},(-\infty,0)\in D_{\varepsilon}, $ where we the
functions are regular in the upper part of the $\mathbb C.$

4.  The function $\mathcal F_0\mathcal F_0(Z (x))$ is regular in
the upper part of $\mathbb C$ for all $p:\{Im\,p>0\}.$ and it is
continuous for all $p=x\in(-\infty,0),$ (if $1<l>n+2.)$ We use
$|e^{ipt}|\leq e^{-y},p=x+iy,y>0.$

\centerline{ References}

1. Pavlov A.V.The random series of Fourier  and its apply to the
theory of prognosis-filtering.(In russian ). Moscow University of
Lomonosov , mechanical-mathemat. faculty, 2000. ISBN
5-93839-002-8, 64 p.

2. Lavrentiev M.A, Shabat B.V. Methods of the theory of functions
of
 complex variables. Moscow :Edd.Science.,1987.-688 p.

\newpage

Home ad.: Prof. Pavlov Andrey
Valerianovich.(Russia,Moscow,109444,Fergankaia st.,11-1,292 )

Post ad.: for correspond.

Russia,Moscow,117454,pr.Vernadskogo,78,MIREA(TU),Cybernet.higher
 mathemat.,\\for prof.Pavlov Andrey.Valerianovich.

E-mail:

AVpavlov@rambler.ru

and11pavlov@msn.com

(095)4330355-job.tel. in Moscow.

%\pagestyle{empty}
%\oddsidemargin=5mm
%\textwidth=17cm
%\topmargin=-17mm
%=-17mm
%\textheight=255mm
%\usepackage[russian]{babel}
%\usepackage[cp866]{input enc}
%\begin{document}

\end{document}